\documentclass[11pt,doublespace]{article}

\usepackage{graphicx}

\usepackage{textcomp}
\usepackage{amsmath}
\usepackage{amsfonts}
\usepackage{epsfig}

\usepackage{amsfonts}
\usepackage{amssymb}

\renewcommand{\baselinestretch} {1.3}

\makeatletter 
\def\singlespace{\def\baselinestretch{1}\@normalsize}

\newcommand{\reals}{\ensuremath{{\mathbb R}}}
\newcommand{\RR}{\reals}

\newcommand{\EE}{\ensuremath{{\mathbb E}}}
\newcommand{\PP}{\ensuremath{{\mathbb P}}}

\newtheorem{theorem}{Theorem}[section]

\newtheorem{lemma}{Lemma}[section]
\newtheorem{defin}{Definition}[section]

\newtheorem{example}{Example}[section]
\newtheorem{remark}{Remark}[section]
\newtheorem{col}{Corollary}[section]

\DeclareGraphicsExtensions{.eps, .jpg}

\title{{\sc Some new approaches to infinite divisibility}\footnote{The authors wish to dedicate this article to Professor Calyampudi Radhakrishna Rao on his 90th birthday, as a mark of their respect and admiration for him. }}

\author{ {\sc Theofanis} ~{\sc Sapatinas}, \; {\sc Damodar N.} ~{\sc Shanbhag} \\
Department of Mathematics and Statistics,
          University of Cyprus,\\
          P.O. Box 20537,
          CY 1678 Nicosia, Cyprus\\
          Email: \texttt{fanis@ucy.ac.cy},~
          Email: \texttt{damodarshanbhag@hotmail.com} 
          \vspace{.2cm}
          \\ and \\
          {\sc Arjun K.} ~{\sc Gupta}\\
          Department of Mathematics and Statistics, Bowling Green State University,\\
          Bowling Green, Ohio 43403-0223, USA\\
          Email: \texttt{gupta@bgsu.edu}
}

\date{}

\begin{document}
\pagestyle{empty}
\maketitle
\begin{abstract}
Using an approach based, amongst other things, on Proposition 1 of Kaluza (1928), Goldie (1967) and, using a different approach based especially on zeros of polynomials, Steutel (1967) have proved that each nondegenerate distribution function (d.f.) $F$ (on $\RR$, the real line), satisfying $F(0-) = 0$ and  $F(x) = F(0) + (1-F(0)) G(x)$, $x > 0$, where $G$ is the d.f. corresponding to a mixture of exponential distributions,
 is infinitely divisible. Indeed, Proposition 1 of Kaluza (1928) implies that any nondegenerate discrete probability distribution  $\{p_x: x= 0,1, \ldots\}$ that is log-convex or, in particular, completely monotone, is compound geometric, and, hence, infinitely divisible. Steutel (1970),  Shanbhag \& Sreehari (1977) and Steutel \& van Harn (2004, Chapter VI) have given certain extensions or variations of one or more of these results.  Following a modified version of the C.R. Rao {\em et al.} (2009, Section 4) approach based on the Wiener-Hopf factorization, we establish some further results of significance to the literature on infinite divisibility.
 
\vspace{2mm}

{\bf Key words and phrases}: {Kaluza sequences; Infinite
divisibility; Log-convexity; Mixtures of exponential distributions;
Mixtures of geometric distributions; Wiener-Hopf factorization}

\vspace{2mm}

{\bf AMS (2000) Subject Classification}: {Primary 60E05; secondary
62E10}

\end{abstract}
\pagestyle{plain}

\newpage

\section{Introduction}
\label{sec:intro}

Several papers, including those of Goldie (1967) and Steutel (1967),
address the preservation of infinite divisibility property under
mixing. Steutel (1970) and Steutel \& van Harn (2004) have reviewed
and unified most of the literature in this respect. This consists of, amongst other things, the result, due to both Goldie (1967) and Steutel (1967), that each nondegenerate d.f. $F$ (on $\RR$, the real line), satisfying $F(0-) = 0$ and $F(x) = F(0) + (1-F(0)) G(x)$, $x > 0$, where $G$ is the d.f. corresponding to a mixture of exponential distributions (or, in other words, corresponding to a density with its restriction to $(0, \infty)$ completely monotone), is infinitely divisible (i.d.), and some of its extensions and variations. It may be worth pointing out here that the result relative to exponential distributions mentioned above is usually referred to as the Goldie-Steutel theorem or Goldie-Steutel result; to be short, we refer henceforth to this simply as {\tt GSR}.

Proposition 1 of Kaluza (1928), a renewal theoretic result, implies clearly that any log-convex or, in particular, completely monotone, nondegenerate discrete probability distribution $\{ p_x: x=0,1,\ldots\}$ is compound geometric, and, hence, i.d. In view of the information in Steutel (1970, p. 89) (and the closure property of the class of i.d. distributions), one can then see that, as a corollary to Kaluza's result, any nondegenerate d.f., $F$ (concentrated) on $[0, \infty)$, that is differentiable, with log-convex derivative, on $(0, \infty)$, is i.d.  (In Section 2 to follow, we are to discuss some basic properties of log-convex sequences and log-convex functions, of relevance to this study.) The two results mentioned here obviously extend {\tt GSR} and its discrete version (respectively, in the reverse order). In his argument to obtain {\tt GSR}, Goldie (1967) also uses the proposition referred to, but, he does so in conjunction with a certain characteristic property of the i.d. distributions on $[0, \infty)$, based, essentially, on mixtures of Poisson distributions. (We view here the degenerate distribution at the origin as Poisson with mean zero.)  Incidentally, the characterization, just referred to, shows especially, that {\tt GSR} implies its discrete version and vice versa (and that this is also so for Theorems VI.3.10 and VI.7.9 of Steutel \& van Harn (2004), respectively).

Shanbhag \& Sreehari (1977, Theorem 5) and Steutel \& van Harn (2004, Theorems VI.3.13 and
VI.7.10) provide us with extended versions and variations of {\tt GSR}; in these contributions, Theorem 2.3.1 of Steutel (1970) has implicitly played an important role or provided a motivation. Also, it may be noted here that Theorem VI.3.13 of Steutel \& van Harn (2004) is a corollary to Theorem 5 of Shanbhag \& Sreehari (1977), and that  some applications of {\tt GSR} and its extensions have appeared especially in Shanbhag \& Sreehari (1977) and Shanbhag {\em et al.} (1977). The main theme of the arguments used by Steutel (1967, 1970) and Steutel \& van Harn (2004) to obtain the results implied is obviously based on zeros of polynomials, but, more recently, C.R. Rao {\em et al.} (2009, Section 4) have, effectively, shown that these results are also by-products of the Wiener-Hopf factorization met in the theory of random walk, thanks to Kaluza (1928, Proposition 1).

Obviously, {\tt GSR} asserts that if $W$ and $X$  are independent random variables (r.v.'s) with $W$ nonnegative and $X$ exponential, then $WX$ is i.d., while each of Theorem 5 of Shanbhag \& Sreehari (1977) and Theorem VI.3.13 of Steutel \& van Harn (2004) asserts that the result in question remains valid irrespectively of whether or not the r.v. $W$ is nonnegative. Moreover, the discrete version of {\tt GSR} is extended by Theorem VI.7.10 of Steutel \& van Harn (2004), which, in turn, states, implicitly, that if $X$ is an integer-valued r.v. with $\PP\{X=0\} > 0$ and the conditional distributions of $X$ given, that $X$ is nonnegative and of $-X$ given that $-X$ is nonnegative, are both completely monotone (on $\{0,1,\ldots\}$), then $X$ is i.d.

In this article, in Section \ref{sec:newsection}, we present two auxiliary lemmas and, in Section \ref{sec:w-h}, we extend the argument based on the Wiener-Hopf factorization, met in C.R. Rao {\em et al.} (2009), to establish certain results based on log-convex sequences and log-convex functions, respectively, and demonstrate, also, that Theorems VI.7.10 and VI.3.10 of Steutel \& van Harn (2004) are their corollaries. Besides, we include, in both of these sections, various remarks with comments, comprising examples in some cases, on our main findings as well as on results in the existing  literature, that are of relevance to the present investigation.

\section{Auxiliary lemmas and related observations}
\label{sec:newsection}
We may start this section by stating first the following crucial definition:

\begin{defin}
\label{def1_FF} Let $B \subseteq \RR$ such that $B=\{a,a+b,a+2b,\ldots\}$ or $B=(a,\infty)$ with $a \in \RR$ and $b \in (0,\infty)$. Then, any function $g: B\mapsto [0,\infty)$, such that
\begin{equation}
\label{eq:log_convF}
g^2(x+y) \leq g(x) \,g(x+2y), \quad \text{for each} \quad  x \in B \quad \text{and} \quad y \in C,
\end{equation}
where 
\begin{equation*}
\label{eq:cont} C= \left\{
\begin{array}{ll}
\{b\} & {\rm if}\;\;\; B=\{a,a+b,a+2b,\ldots\}, \\
(0,1) & {\rm if}\;\;\;
B=(a,\infty),
\end{array} \right.
\end{equation*}
will be referred to as a log-convex function (on $B$).
\end{defin}

\noindent (Note that, for notational convenience, in Definition \ref{def1_FF}, we allow $g$, that is identically equal to 0, to be also called log-convex, and that Zygmund (2002, p. 25) has used the same convention in his definition of a log-convex function.) 

Obviously, by Definition \ref{def1_FF}, a log-convex function $g$ that is not identically equal to 0 on $B\setminus\{a\}$ is non-vanishing, and, essentially, by (an observation in) Lo$\grave{{\rm e}}$ve (1963, p. 159), in the case of $B=(a,\infty)$, a log-convex function (on $B$), as defined in Definition \ref{def1_FF}, is either continuous or not a Borel function; for the literature supporting the claim of Lo$\grave{{\rm e}}$ve (1963, p. 159), one may consult, e.g., Hardy {\em et al.} (1952, p. 96) or Donoghue (1969, Section 1.3).  It may also be pointed out that, if $g$ is a log-convex function, then given any $x$ and $y$ as in (\ref{eq:log_convF}), we get the restriction to $\{x,x+y,x+2y,\ldots\}$ as log-convex on that set, and we refer to the sequence $\{g(x+ny):\, n=0,1,\ldots\}$ as a log-convex sequence.

\medskip

In what follows, in this section, we give the auxiliary lemmas and remarks  of relevance to these.

\begin{lemma}
\label{lem1_FN}
Let $B$ be as in Definition \ref{def1_FF}, let $(\Omega, \cal{E}, \mu)$ be a measure space, and let $\{g_x:\, x \in B\}$ be a family of integrable functions on the measure space so that, for almost all (for short, a.a.) $\omega \in \Omega$, $g_x(\omega)$, $x\in B$, are log-convex. Define (in the notation of Lo$\grave{{\rm e}}$ve (1963, p. 119)), 
$$
h(x) = \int g_x, \quad x \in B.
$$
Then, $h$ is a log-convex function (on $B$). (We do not impose here the restriction that $h$ be continuous if $B=(a,\infty)$.)
\end{lemma}

\noindent {\bf Proof.} If $\alpha, \beta, \gamma \in (0,\infty)$, then, for instance, by the first statement in the last paragraph of Kingman (1972, p. 18), we have 
\begin{equation}
\label{eqeq_FN}
\beta^2 \leq \alpha \gamma
\end{equation}
 to be equivalent to that
\begin{equation}
\label{eq:equivlog_con}
\alpha \lambda^2 - 2 \beta \lambda + \gamma \geq 0, \quad \text{for all} \quad \lambda >0.
\end{equation}
In view of the assumptions of the lemma and the criterion of (\ref{eq:log_convF}) for a function to be log-convex, it follows that if $x \in B$ and $y \in C$ (as in (\ref{eq:log_convF})), we have (\ref{eq:equivlog_con}) to be valid with $\alpha$, $\beta$ and $\gamma$ replaced, respectively, by 
$g_x(\omega)$, $g_{x+y}(\omega)$ and $g_{x+2y}(\omega)$ for a.a. $\omega \in \Omega$, and, hence, to be valid with $\alpha$, $\beta$ and $\gamma$ replaced, respectively, by $h(x)$, $h(x+y)$ and $h(x+2y)$. This, in turn, shows, since (\ref{eq:equivlog_con}) implies (\ref{eqeq_FN}), that (\ref{eq:log_convF}) holds with $h$ in place of $g$, and, consequently, that the lemma holds. \hfill $\Box$

\begin{lemma}
\label{lem2_FN}
Let $G$ be a Lebesgue-Stieltjes measure function (i.e., a nondecreasing, right-continuous, real-valued function) on $\RR$, so that, for some $a \in \RR$, $G(x)=0$ if $x<a$ and $G(x)$ is differentiable with log-convex derivative if $x >a$. Define the sequence $\{G_n:\, n=1,2,\ldots\}$ of Lebesgue-Stieltjes measure functions on $\RR$, so that, for each $m,n \in \{1,2,\ldots\}$,
\begin{equation*}
\label{eq:cont} G_n(x)= \left\{
\begin{array}{ll}
G(a+m/n) & {\rm if}\;\;\; x \in [a+(m-1)/n, a+m/n), \\
0 & {\rm if}\;\;\;
x \in (-\infty,a).
\end{array} \right.
\end{equation*}
Then, $\{G_n:\, n=1,2,\ldots\}$ converges weakly to $G$ and, for each $n \in \{1,2,\ldots\}$, $G_n$ is concentrated on $B_n$ with $G_n(x)-G_n(x-)$, $x \in B_n$, log-convex, where $B_n=\{a,a+1/n,a+2/n,\ldots\}$.
\end{lemma}

\noindent {\bf Proof.} That $\{G_n:\, n=1,2,\ldots\}$ converges weakly to $G$ is obvious, since, by assumptions, for each $n \in \{1,2,\ldots\}$, 
\begin{equation*}
\label{eq:cont} |G_n(x)-G(x)|=G_n(x)-G(x) \left\{
\begin{array}{ll}
\leq G(x+1/n) -G(x) & {\rm if}\;\;\; x \in [a, \infty), \\
=0 & {\rm if}\;\;\;
x \in (-\infty,a).
\end{array} \right.
\end{equation*}
It is also obvious, by assumptions, that, for each $n \in \{1,2,\ldots\}$, $G_n$ is concentrated on $B_n$ and, in view of the observation below 11.82 in Titchmarsh (1978, p. 368), that
\begin{equation}
\label{eqLS_FN}
G_n(x)-G_n(x-)=G(a)S(x)+\int_0^{1/n}G'(x+y)\,dy \quad \text{if} \quad x \in B_n,
\end{equation}
where $S(a)=1$ and $S(x)=0$ if $x>a$. Clearly, we have in (\ref{eqLS_FN}), the function $G(a)S(x)$, $x \in B_n$, to be log-convex, and, in view of the log-convexity of $G'$ on $(a,\infty)$, for each $y \in (0,\infty)$, the function $G'(x+y)$, $x \in B_n$, to be log-convex. Consequently, by Lemma \ref{lem1_FN}, it follows that $G_n(x)-G_n(x-)$, $x \in B_n$, is log-convex, and, thus, the lemma holds. \hfill $\Box$

\begin{remark}
\label{rem2_1F}
{\rm In the literature, usually, a log-convex sequence $\{u_n\}$, with $u_0=1$ and $0 <u_n \leq 1$ if $n >0$, is called a Kaluza sequence, especially, in recognition of the findings on such sequences given in  Kaluza (1928); to be short, we refer to this as a {\tt KS}. Indeed, Proposition 1 of Kaluza (1928) establishes that each {\tt KS} is renewal. Kingman (1972, Section 1.5), in particular, and Shanbhag (1977) amongst others, have made some further observations on these sequences. [Incidentally, the statement in Kingman (1972, p. 18) that the class of {\tt KS}'s is a closed convex subset of the class of renewal sequences requires the convention of the sequence $\{u_n\}$, with $u_0=1$ and $u_n =0$ if $n >0$, being Kaluza, adopted. A similar blemish, also, exists in the {\tt KS}-related proof, given by Goldie (1967) for his Theorem 2, since it assumes, implicitly (in the notation used in it) that $\PP\{Z=0\}<1$, without stating that this is so.] If $g$ is a non-vanishing log-convex function on $B$, then (\ref{eq:log_convF})  implies that, for each pair $x$ and $y$ as in (\ref{eq:log_convF}), the sequence
$\{ g(x+ny)/g(x+(n-1)y):\, n=1,2,\ldots\}$,
is nondecreasing, and, hence, for $n \in \{1,2,\ldots\}$,
$$
\frac{g(x+ny)}{g(x)} = \prod_{m=1}^n \frac{g(x+my)}{g(x+(m-1)y)}
\leq \prod_{m=n+1}^{2n} \frac{g(x+my)}{g(x+(m-1)y)}= \frac{g(x+2ny)}{g(x+ny)},
$$
implying that
\begin{equation}
\label{FF_llcc}
g^2(x+ny) \leq g(x)\,g(x+2ny).
\end{equation}
In view of (\ref{FF_llcc}), it follows that (\ref{eq:log_convF}) is equivalent to its version with $C$ replaced by $C^*$, where $C^*=\{b,2b,\ldots\}$ in the discrete case, and $C^*=(0,\infty)$ otherwise.
Obviously, there are several other equivalent formulations for (\ref{eq:log_convF}).}
\end{remark}

\begin{remark}
\label{rem3_2F}
{\rm In view of the result of Titchmarsh (1978) met in the proof of Lemma \ref{lem2_FN}, it follows, by Lemma \ref{lem1_FN}, that $\{ n ( G (x+(1/n)) - G (x)), x>a : n =1,2,\ldots\}$, with $G$ as in Lemma \ref{lem2_FN}, is a sequence of continuous log-convex functions; since the sequence converges (pointwise) to the function $G' (x )$ , $x > a$, it is seen, by (10.4) and (10.7) in Volume I, Chapter 1, of Zygmund (2002), that this latter function is indeed a continuous log-convex function. That $G'$ is continuous,  follows, also, from Lo$\grave{{\rm e}}$ve (1963, p. 159), since it is a log-convex Borel function; we could have, obviously, used this information, in place of that in Titchmarsh (1978), in the proof of Lemma \ref{lem2_FN}. It may be noted, in this connection, that Lo$\grave{{\rm e}}$ve (1963) and Zygmund (2002) use different approaches to define convex functions, though, these turn out, in view of (10.7) of Zygmund (2002), to be equivalent in the case of continuous functions defined on open intervals, and we follow (in Definition \ref{def1_FF}) implicitly Lo$\grave{{\rm e}}$ve (1963) to define a log-convex function $g$ in the case of $B = (a,\infty)$.}
\end{remark} 

\begin{remark}
\label{rem2_2F}
{\rm $(i)$ If $B=\{a,a+b,a+2b,\ldots\}$ with $a \in \RR$ and $b \in (0,\infty)$, then, by Lemma \ref{lem1_FN}, any nonnegative real function $h$ on $B$, of the form
\begin{equation}
\label{eq:rem22F}
h(x) =\sum_{n=0}^{\infty} g(x+nb) \nu_n, \quad x \in B,
\end{equation}
where $g$ is a log-convex function on $B$ and $\{\nu_n:\, n=0,1,2,\ldots\}$ is a sequence of nonnegative reals, is log-convex. 

$(ii)$ Suppose $B = (a, \infty)$ with $a \in \RR$, and $g$ is a continuous log-convex function on $B$. Then, in view of Lemma \ref{lem1_FN}, by Lo$\grave{{\rm e}}$ve (1963, p. 159), it follows that (since it is a Borel function) any nonnegative real function on $B$, satisfying
\begin{equation}
\label{eq:rem2iiFFF}
h(x) = \int_{[0,\infty)}g(x + y) d\nu(y),\quad  x \in B,                                                   
\end{equation}
with $\nu$ as a Lebesgue-Stieltjes measure, is a continuous log-convex function. Similarly, if $a$ nonnegative, using the relevant information (i.e., relative to $C^*$) appearing in Remark \ref{rem2_1F}, it can be seen that any nonnegative real function $h$ on $B$, satisfying (\ref{eq:rem2iiFFF}), with ``$x +y$'' replaced  by ``$x y$'' and ``$[0,\infty)$'' replaced, respectively, by ``$(0, \infty)$'' if $a=0$ and ``$[1,\infty)$'' otherwise,  is also a continuous log-convex function. (Since, in each case, $\{h_n:\, n=1,2,\ldots\}$, where, for each $n$, $h_n$ satisfies the relevant version of (\ref{eq:rem2iiFFF}), with $\nu (\cdot\cap [0,n])$ in place of $\nu$, is a sequence of continuous functions converging, by the monotone convergence theorem, to $h$, that $h$ has the stated property is implied, in view of Lemma \ref{lem1_FN} and the last observation in Remark \ref{rem3_2F}, also, by Zygmund (2002; Theorem 10.4).)

$(iii)$ If $a$ is nonnegative and $H$ is a Lebesgue- Stieltjes measure function on $\RR$, such that
\begin{equation}
\label{eq:rem2iiF}   
H(x)= \left\{
\begin{array}{ll}
\EE ( G ( x V)) & {\rm if}\;\;\; x \geq a, \\
0 & {\rm otherwise,}
\end{array} \right.
\end{equation}                                                        
where $G$ is as in Lemma \ref{lem2_FN} and $V$ is a nonnegative real r.v. meeting, if $a>0$, a further condition of $\PP\{V \geq 1\}=1$, it follows (on reading, for convenience, $0\, G'(0)=0$), by Fubini's theorem, that 
\begin{equation}
\label{eq:rem2iiFF} 
H(x) - H(a) =  \int_{(a,x)} \EE ( V G' (y V) ) dy , \quad   x \in (a, \infty).    
\end{equation}                 
Clearly, then, in view of the properties of $G'$ referred to in Remark \ref{rem3_2F} and our observation on the latter version of $h$ in $(ii)$ above, we have that $H$ is differentiable on $(a, \infty)$ with continuous log-convex derivative ($E( V G'(xV))$ for $x> a$).  
}
\end{remark}

\begin{remark}
\label{rem:2F} {\rm Proposition 1 of Kaluza (1928) obviously implies that any nondegenerate log-convex 
probability function  on  $\{a, a+b, a+2b,\ldots\}$, with $a \in \RR$ and $b \in (0, \infty)$, corresponds to an  r.v. $a +b X$, where $X$ is an r.v. with compound geometric distribution, and, hence, is i.d. In view of this, Lemma \ref{lem2_FN} implies that if $G$ in it is a d.f. or $H$ of Remark \ref{rem2_2F} $(iii)$ is a d.f., then it is i.d.; this latter result is an extension of Theorem 4.2.6 of Steutel (1970) (or of Theorem III. 8.4 in Steutel and van Harn (2004)). The two results based on Kaluza's proposition, met here, clearly extend (respectively, in the reverse order)  GSR and its discrete analogue. We have already given some information to this effect in the introduction. In view of this, C.R. Rao {\em et al.} (2009, Section 4), tells us that Theorem VI.3.13 and Theorem VI.7.10 of Steutel \& van Harn (2004) have alternative proofs  based, at least partly, on {\tt KS}'s or, in particular, completely monotone sequences. 
}
\end{remark}

\section{Results based on log-convex sequences and log-convex functions}
\label{sec:w-h}

To recall partially the information already implied, in our discussion so far, especially in Remark \ref{rem:2F}, applying the  Wiener-Hopf factorization relative to a random walk,  C.R. Rao {\em et al.}  (2009, Section 4) effectively proved Theorems VI.3.13 and VI.7.10 of Steutel \& van Harn (2004). As hinted in the remark referred to, this approach relies also upon certain properties of log-convex sequences, or, in particular, of  {\tt KS}'s. Included in the following two subsections, viz., Subsections \ref{subsec:3.1} and \ref{subsec:3.2}, of the present section, are some further {\tt KS}-related results and remarks. In Subsection \ref{subsec:3.1}, modifying the relevant arguments in C.R. Rao {\em et al.} (2009) appropriately, we establish a key theorem, viz., Theorem \ref{th:w-h}, with assertion in it based on log-convex sequences, and, then, give four of its important corollaries and, in Subsection \ref{subsec:3.2}, make some pertinent observations through remarks on the results obtained. From the presentation of Subsection \ref{subsec:3.1}, it is clear that Corollaries \ref{col:1} and \ref{col:3}, with assertions based on log-convex sequences and log-convex functions, respectively, imply Corollaries \ref{col:2} and \ref{col:4}, respectively, and, also, that the latter corollaries are respective rephrased versions of the  aforementioned Steutel \& van Harn (2004) theorems.
 
\subsection{The key theorem and its corollaries}
\label{subsec:3.1}
 
We now present in this subsection our main result, viz., Theorem \ref{th:w-h} and four of its corollaries, referred to above.

\begin{theorem}
\label{th:w-h} Let $F$ be a d.f. relative to a
probability distribution $\{p_x:~x=0,\pm1,\pm2,\ldots\}$
such that, for some constant $K>0$,
\begin{equation}
\label{eq:disc} K p_x = \left\{
\begin{array}{ll}
\sum_{j=|x|+1}^{\infty}  v_{-j}  & {\rm if}\;\;\; x =-1,-2,\ldots, \\
\sum_{j=x+1}^{\infty}  v_j & {\rm if}\;\;\;
 x=1,2,\ldots, \\
\max \big\{\sum_{j=1}^{\infty}  v_{-j}, \sum_{j=1}^{\infty}  v_j \big\} 
& {\rm if}\;\;\; x=0,
\end{array} \right.
\end{equation}
with $\{v_j:~j=\pm 1, \pm2,\ldots\}$ as a sequence of nonnegative reals
for which  $\{v_{-j}:~j=1,2,\ldots\}$ and $\{v_j:~j=1, 2,\ldots\}$
are log-convex. Then, $F$ is i.d.
\end{theorem}

\noindent {\bf Proof.} We may assume, without loss of generality,
that $\{v_j:~j=0,\pm1,\pm2,\ldots\}$, with $v_0$ denoting some nonnegative real number, is so that $\sum_{j=1}^{\infty}  v_{-j}=\sum_{j=1}^{\infty}v_j$ and
$\sum_{j=-1}^{1}v_j=1$; this follows since if $\sum_{j=1}^{\infty}v_{-j} \neq \sum_{j=1}^{\infty}v_j$, then, taking, without loss of generality, $\sum_{j=1}^{\infty}v_{-j}>\sum_{j=1}^{\infty}v_{j}$, we can 
replace $v_1$ by $v_1+\sum_{j=1}^{\infty}v_{-j}-\sum_{j=1}^{\infty}v_{j}$ and verify that the normalized version of the resulting $v$-sequence is as required and
(\ref{eq:disc}) (possibly, with different $K$) holds with the original $v$-sequence replaced by this.
Now, since we need a proof for the theorem only when the distribution (satisfying (\ref{eq:disc})) is  nondegenerate, we can assume, again, without loss of generality, that $v_j  > 0$ for all $j>0$, and
define, for each positive integer $k$, a sequence
$\{v_j^{(k)}:~j=0,\pm1,\pm2,\ldots\}$ satisfying
\begin{equation*}
\label{eq:w-h-new} v_j^{(k)} = \left\{
\begin{array}{ll}
v_j & {\rm if}\;\;\; j=0,-1,\pm 2,\ldots,\pm k, -(k+1), -(k+2),\ldots, \\
v_k \big( \frac{v_{k+1}}{{v_k}} \big)^{j-k} & {\rm if}\;\;\; j =k+1,k+2,\ldots,\\
v_1 + \sum_{i=2}^{\infty} \big(v_i -
v_i^{(k)} \big) & {\rm if}\;\;\; j =1.
\end{array} \right.
\end{equation*}

Kingman (1972, p. 18) involves the idea that we have used in the construction of $\{v_j^{(k)}\}$ from
$\{v_j\}$ and implies, in view of (\ref{eq:disc}), that $0<v_{k+1}/v_k <1$; note, in particular, that, in the present case, 
$\{v_{1+j}/v_1: j=0,1,2,\ldots\}$ and $\{v_{-1-j}/v_{-1}: j=0,1,2,\ldots\}$ are such that the first one is a decreasing {\tt KS} and, unless it is the sequence $\{1,0,0,\ldots\}$, so also is the second one. Obviously, for each positive integer $k$, $\{(v_{-1}^{(k)}+v_{0}^{(k)}+v_{1}^{(k)})^{-1} v_{j}^{(k)}: j=0, \pm 1, \pm 2,\ldots\}$ is a sequence possessing 
the properties sought of $\{v_j: j=0, \pm 1, \pm 2, \ldots\}$ above and provides us with a distribution $\{p_x^{(k)}\}$ satisfying the version of (\ref{eq:disc}) (possibly, with different $K$ and) with $v_j$'s replaced by $v_j^{(k)}$'s. Denoting the d.f. relative to this latter distribution by $F_k$, it is easily seen that $\{F_k: k=1,2,\ldots\}$ converges weakly to $F$, and hence, in view of the closure property of the class of
i.d. distributions under weak convergence, it is sufficient if we prove that, for each positive integer $k$,  $F_k$ is i.d.

Consequently, it follows that to prove that $F$ is i.d., there is no loss of generality if we assume that for some $j_0 \in \{1,2,\ldots\}$
and $b \in (0,1)$, $v_j \propto b^j$, $j \geq j_0$. Assume then that $\{v_j: j=0,\pm 1, \pm 2,\ldots\}$ meets this additional condition and define $c=-\ln b$. Also, define now the sequence $\{w_j: j=0,\pm 1, \pm 2, \ldots\}$ so that
\begin{equation}
\label{eq:w-new-FF} w_j = \left\{
\begin{array}{ll}
v_j-v_{j-1} & {\rm if}\;\;\; j =-1,-2,\ldots,\\
v_j-v_{j+1} & {\rm if}\;\;\; j=1, 2,\ldots, \\
v_0 & {\rm if}\;\;\; j =0,
\end{array} \right.
\end{equation}
and observe that this is a distribution with mean 0, having an m.g.f. with its domain of definition as a superset of $[0, c)$.
Obviously, in the present case, $\{p_x\}$ has the same property as $\{w_j\}$ that it has an m.g.f. with $[0,c)$ as a subset of its domain of definition; denote then the m.g.f.'s of $\{p_x\}$ and $\{w_j\}$ by $M$ and $M^*$, respectively. Following the relevant approach of C.R. Rao {\em et al.} (2009, Section 4.2) (involving Theorem XII.2.2 of Feller (1971) or otherwise), we can now conclude that the weak descending and the ascending ladder height measures associated with the random walk relative to $\{w_j\}$ are indeed probability measures with m.g.f.'s having domains of definitions as supersets of $[0,c)$, and that 
\begin{equation}
\label{eq:fanisNew}
K M(t) =  \frac{(1-M^*(t))}{(1-e^{-t})(1-e^{t})} = 
\left(\frac{1-M^*_1(t)}{1-e^{-t}}\right)\left(\frac{1-M^*_2(t)}{1-e^{t}}\right),
\quad t \in (0, c),
\end{equation}
where $M^*_1$ and $M^*_2$ are, respectively, the m.g.f.'s of the weak descending and the ascending ladder height
measures referred to. (For a simple argument to see that the first equation in (\ref{eq:fanisNew}) holds, refer to Remark \ref{rem:2} $(i)$, given below.) In view of (XII.3.7a) in Feller (1971) corresponding to the distribution of a weak descending ladder height and its analogue corresponding to an ascending ladder height, it follows that, in our case, the distributions referred to are of the forms $\{\sum_{n=0}^{\infty}w_{i-n} \nu^{(1)}_n : i = 0,-1,-2,\ldots\}$ and $\{ \sum_{n=0}^{\infty} w_{i+n} \nu^{(2)}_n: i=1,2,\ldots\}$, with $\nu^{(1)}_n$ and $\nu^{(2)}_n$ nonnegative for all $n$, respectively. Since for any nonnegative integer-valued random variable $Z$, we have a standard result that, for each $t \neq 0$ for which $\EE(e^{tZ}) < \infty$, 
$$
(1- \EE(e^{t Z}))/(1-e^t) = \sum_{j=0}^{\infty}  e^{t j} \PP\{Z >j\},
$$
on appealing to (\ref{eq:fanisNew}), especially, in view of Fubini's theorem, Remark \ref{rem2_2F} $(i)$ and the log-convexity properties of
$\{v_{-j}:~j=1,2,\ldots\}$ and $\{v_{j}:~j=1,2,\ldots\}$, it is clear that
there exist sequences $\{v_{1j}^*:~j=0,1,2,\ldots\}$ and
$\{v_{2j}^*:~j=0,1,2,\ldots\}$ of nonnegative reals and of positive reals, respectively, that are log-convex
such that
\begin{equation}
\label{eq:fanisNewFF}
K M(t) = \Bigg( \sum_{j=0}^{\infty}v_{1j}^*e^{-tj}\Bigg)
\Bigg( \sum_{j=0}^{\infty}v_{2j}^*e^{tj}\Bigg), \quad t \in [0,c).
\end{equation}
(Note that we have assumed in this proof, without loss of generality, that $v_j >0$ for all $j>0$.) We have, obviously, $\{v_{1j}^*:~j=0,1,2,\ldots\}$ and
$\{v_{2j}^*:~j=0,1,2,\ldots\}$ in (\ref{eq:fanisNewFF}) to be, in view of the related information in Remark \ref{rem2_1F}, proportional to renewal sequences and, hence, to appropriate discrete i.d. distributions (with the first one as degenerate at the origin or compound geometric, and the second one as compound geometric). We have then that $F$ is i.d., implying that
the theorem holds. \hfill $\Box$

\begin{col}
\label{col:1}
For each $\alpha \in [0,1]$,  let  $F^{(\alpha)}$ denote a d.f. so that
\begin{equation}
\label{eq:fanalpha}
F^{(\alpha)} (x)= \alpha F_1(x) + (1-\alpha) F_2(x), \quad x \in \RR,
\end{equation}
where $F_1$ is the d.f. relative to the degenerate distribution at the origin and $F_2$ 
is the d.f. of a nondegenerate distribution satisfying (\ref{eq:disc}). Then,  each $F^{(\alpha)}$ is i.d.
\end{col}

\noindent {\bf Proof.} We need a proof only when $\alpha \in (0,1)$. Note now that, in this case, for each $\alpha$,
the probability function (on $\{0, \pm 1, \pm 2, \ldots\}$) relative to $F^{(\alpha)}$ satisfies (\ref{eq:disc}),  with $K$ unaltered and $\{v_j\}$ replaced by $\{v_j^{(\alpha)}\}$, where

\begin{equation*}
\label{eq:fanisCor1}
v_j^{(\alpha)} = \left\{
\begin{array}{ll}
K \alpha+(1-\alpha)v_j & {\rm if}\;\;\; j=-1, 1, \\
(1-\alpha)v_j &  {\rm otherwise}.
\end{array} \right.
\end{equation*}
Consequently, we have the corollary. \hfill $\Box$

\begin{col}
\label{col:2}
Each probability distribution $\{p_x: x=0, \pm 1, \pm 2, \ldots\}$, with $p_0 >0$ and 
$\{p_{-x}: x=0, 1, 2, \ldots\}$ and $\{p_x: x=0, 1, 2, \ldots\}$  completely monotone, is i.d.
\end{col}

\noindent {\bf Proof.} In view of Hausdorff's theorem, referred to, e.g., in Theorem VII.3.1 of Feller (1971), it  follows that
\begin{equation*}
\label{eq:fanisCor2}
p_x = \left\{
\begin{array}{ll}
p_0 m_{|x|}^{(1)} & {\rm if}\;\;\; x=-1,-2,\ldots,\\
p_0 m_x^{(2)} & {\rm if}\;\;\; x=1,2,\ldots, \\
p_0 m_0^{(1)} (=p_0 m_0^{(2)}) &  {\rm if}\;\;\; x=0, 
\end{array} \right.
\end{equation*}
where $\{m_x^{(1)}: x=0,1,2,\ldots\}$ and $\{m_x^{(2)}: x=0,1,2,\ldots\}$ (with, obviously, $m_0^{(1)}=m_0^{(2)}=1$) are moment sequences
relative to probability distributions concentrated on $[0,1)$. Defining  $\{v_j: j=\pm 1, \pm 2, \ldots\}$ for which
\begin{equation*}
\label{eq:fanisCor2:2}
v_{j+1} = \left\{
\begin{array}{ll}
p_{j+2} -p_{j+1} &  {\rm if}\;\;\; j=-2,-3,-4,\ldots,\\
p_j -p_{j+1}& {\rm if}\;\;\; j=0,1,2,\ldots, 
\end{array} \right.
\end{equation*}
it is then seen that $\{p_x\}$ satisfies (\ref{eq:disc}) with $v_j$'s meeting the required condition and $K=1$. Hence, we have the corollary. \hfill $\Box$

\begin{col}
\label{col:3}
For each $\alpha \in [0,1]$,  let $F^{(\alpha)}$ be so that (\ref{eq:fanalpha}) is met but for that  $F_2$, now, instead of that in the statement of Corollary \ref{col:1}, is the d.f. relative to an absolutely continuous distribution with density $f_2$ satisfying 
\begin{equation}
\label{eq:fanisCor3}
f_2(x) = \left\{
\begin{array}{ll}
\int_{|x|}^{\infty} v_1(y) dy & {\rm if}\;\;\; x<0,\\
\int_{x}^{\infty} v_2(y) dy & {\rm if}\;\;\; x>0, 
\end{array} \right.
\end{equation}
with $v_1$ and $v_2$ as (nonnegative) log-convex functions on $(0, \infty)$. Then, each $F^{(\alpha)}$ is i.d.
\end{col}

\noindent {\bf Proof.} For each $n \in \{1,2,\ldots\}$, define the d.f. $F_2^{(n)}$ such that it is concentrated on $\{0, \pm\frac{1}{n}, \pm \frac{2}{n}, \ldots\}$, having $\{F_2^{(n)}((x+1)/n)-F_2^{(n)}(x/n):  x=0, \pm 1, \pm 2, \ldots\}$ to be of the form of $\{p_x\}$ of Theorem  \ref{th:w-h}, satisfying 
(\ref{eq:disc}) with $K_n >0$ (where $K_n \rightarrow 1$ as $n \rightarrow \infty$) in place of $K>0$ and $\{v_j\}$ replaced by $\{v_j^{(n)}\}$, where
\begin{equation}
\label{eq:fanisCor3:3}
v_j^{(n)} = \left\{
\begin{array}{ll}
\int_0^{1/n} \int_{0}^{1/n} v_1((|j+1|/n)+y+z)dy dz &  {\rm if}\;\;\; j=-1,-2,\ldots,\\
\int_0^{1/n} \int_{0}^{1/n} v_2(((j-1)/n)+y+z)dy dz & {\rm if}\;\;\; j=1,2,\ldots~.
\end{array} \right.
\end{equation}
Note that, for each $n$,  given $y,z, \in 0,1/n)$,  $\{ v_r(( m/n) +y+z): m=0,1,\ldots\}$, $r=1,2$, are log-convex sequences. Hence, by Lemma \ref{lem1_FN}, it follows that, for each $n$, $\{v_j^{(n)}\}$ meets the requirements of the theorem. Since $F_2^{(n)}$ can easily be seen to be so that it converges to $F_2$ weakly, essentially, in view of Corollary \ref{col:1} and the closure property of the class of i.d. distributions, we can conclude that, for each $\alpha \in [0,1]$, $F^{(\alpha)}$ is i.d., and, thus, we have the corollary. \hfill $\Box$

\begin{col}
\label{col:4}
For each $\alpha \in [0,1]$,  let $F^{(\alpha)}$ be so that (\ref{eq:fanalpha})  is met but for that $F_2$, now, instead of that in the statement of Corollary \ref{col:1}, is the d.f. relative to an absolutely continuous distribution with density $f_2$ for which, for each $r\in\{1,2\}$, the function $f_2((-1)^{r}t)$, $t \in (0,\infty)$, is completely monotone. Then, each $F^{(\alpha)}$ is i.d.
\end{col}

\noindent {\bf Proof.} In view of the conditions to be met by $f_2$, essentially, appealing to a version of Bernstein's theorem appearing as Theorem XIII.4.1a of Feller (1971), it follows that, in this case, (\ref{eq:fanalpha}) holds for each $\alpha \in [0,1]$, with $v_r$'s in (\ref{eq:fanisCor3}) so that
$$
v_r(y)=\int_{(0,\infty)} \lambda^2 e^{-\lambda y} d\mu_r(\lambda), \quad y \in (0,\infty), \quad r=1,2,
$$
where $\mu_r$, $r=1,2$, are measures so that $\mu_1+\mu_2$ is a probability measure concentrated on $(0,\infty)$. Clearly, $v_1$ and $v_2$ in this case are log-convex on $(0, \infty)$, satisfying (in obvious notation) $v_r (.)v_r^{''} (.) - (v_r^{'})^2(.) \geq 0$, $r=1,2$.
Hence, we have the corollary from Corollary \ref{col:3}. \hfill $\Box$

\subsection{Some relevant remarks}
\label{subsec:3.2}

We devote the present subsection, as implied before, to making specific observations on our findings through remarks:

\begin{remark}
\label{rem:1}
{\rm The equation (\ref{eq:disc}) in the statement of Theorem \ref{th:w-h} can be also expressed, with a minor notational adjustment, so as to have $K=1$. However, to have the options such as that in which  $\{v_{-1-j}: j=0,1,2,\ldots\}$ and $\{v_{1+j}: j=0,1,2,\ldots\}$ are moment sequences open to us, and to avoid unnecessary notational complications in the arguments used to prove Theorem \ref{th:w-h}, we have decided to retain $K$ in the equation referred to. Also, in view, especially, of the relevant information in Remarks \ref{rem2_2F} $(i)$, \ref{rem2_2F} $(ii)$ and \ref{rem:2F}, it is obvious that if $v_j =0$ for  some (and hence all) $ j >1$ or $j <-1$, Theorem \ref{th:w-h} and Corollary \ref{col:1}, and, if  $g_1(y)$ or $g_2(y) = 0$ for some (and hence all) $y >0$, Corollary \ref{col:3}, follow as simple consequences of Proposition 1 of Kaluza (1928). To illustrate that Corollaries \ref{col:2} and \ref{col:4} are more restrictive versions of Corollaries \ref{col:1} and \ref{col:3}, respectively, we now give the following simple example:

\begin{example}
\label{ex:0F}
{\rm Let $g$ be a function defined on $(0,\infty)$, such that
\begin{equation}
\label{eq:0F}
g(x) =  e^{- x + h(x)},  \quad x>0,                                                      
\end{equation}
where $h(x)=(1-x)^2$ if $x \in (0,1)$ and $h(x)=0$ if $x \geq 1$. Then, $g$ and $\int_x^{\infty}g(y)dy$, $x>0$, are log-convex but not completely monotone; the log-convexity of the first function is obvious and of the second function follows by the relevant information in Remark \ref{rem2_2F} $(ii)$, and that the functions are not completely monotone follows, since these are not differentiable twice and thrice, respectively, at the point $x=1$. Also, if we now define  $\{v_n:\, n=0,1,2,\ldots\}$ such that 
$v_0=e$ and, for each $n>0$, $v_n=g(n/2)$, where $g$ is as in (\ref{eq:0F}), then, by the Hausdorff theorem, neither of the sequences $\{ v_n\}$  and $\{\sum_{m=n+1}^{\infty} v_{m-1} : n=0,1,2\ldots\}$ is completely monotone; note that the cited theorem implies that a real sequence $\{v^*_n: n=0,1,2,\ldots\}$, for which $e\, v^*_{n+2} -2 e^{1/2}v^*_{n+1}+ v^*_n = 0$ for some $n>1$, turns out to be completely monotone, only if $e\, v^*_3 -2 e^{1/2}v^*_2 + v^*_1= 0$ (and this criterion is not met for the two sequences). However, it is easily seen now that these sequences are indeed log-convex.
}
\end{example}
}
\end{remark}

\begin{remark}
\label{rem:2}
{\rm  $(i)$ A simple argument to show that the first equation in (\ref{eq:fanisNew}) is valid, is as follows:
Clearly, in view of the assumptions in the statement of Theorem \ref{th:w-h} in conjunction with those that we added in the proof of the theorem claiming that there was no loss of generality in doing so we have from
(\ref{eq:disc}) (in obvious notation)
\begin{eqnarray*}
K(1-e^{-t})(1-e^t)M(t)&=& K(2-e^{t}-e^{-t})M(t)\\
&=& \sum_{x=-\infty}^{\infty} e^{tx} [K (2 p_x-p_{x-1}-p_{x+1})]\\
&=& 1 - \sum_{x=-\infty}^{\infty} e^{tx} w_x
=(1-M^*(t)), \quad t \in [0, c).
\end{eqnarray*}
This is obvious on noting, in particular, that the last but one equation holds since 
$$
K(2 p_x-p_{x-1}-p_{x+1})=\left\{
\begin{array}{ll}
v_{x-1}-v_{x}=-w_x & {\rm if}\;\;\; x =-1,-2,\ldots,\\
v_{x+1}-v_{x}=-w_x & {\rm if}\;\;\; x=1, 2,\ldots, \\
v_{-1}+v_1=1-v_0=1-w_0 & {\rm if}\;\;\; x =0.
\end{array} \right.
$$

$(ii)$ The argument that we have used in the proof of Theorem \ref{th:w-h} remains valid with $v_0=0$. However, to make the link between our approach in this case and that appearing in Section 4.2 of C.R. Rao {\em et al.} (2009) more transparent, we have allowed here also the case $v_0 \neq 0$. To illustrate this, we now consider the following example.

\begin{example}
\label{ex:1}
{\rm Let $n \in \{1,2,\ldots,\}$ and $\{p_x: x =0,\pm1,\pm2,\ldots\}$ be a probability distribution, discussed implicitly in Section 4.2 of C.R. Rao {\em et al.} (2009), so that
\begin{equation*}
\label{eq:fanisEx1a}
p_x = \left\{
\begin{array}{ll}
\sum_{s=1}^n c_{1s}p_{1s}^{-x} & {\rm if}\;\;\; x=-1,-2,\ldots, \\
\sum_{s=1}^n c_{2s}p_{2s}^{x} &  {\rm if}\;\;\; x=1,2,\ldots, \\
\max\{c_1^*,c_2^*\} & {\rm if}\;\;\; x=0,
\end{array} \right.
\end{equation*}
where $c_r^*=\sum_{s=1}^n c_{rs}$, $r=1,2$, $c_{rs}>0$ and $p_{rs} \in (0,1)$ for each $r \in \{1,2\}$ and $s \in \{1,2,\ldots,n\}$. Then, by Corollary \ref{col:2}, it is immediate that $\{p_x\}$ is i.d.
(Note that, in this case, the m.g.f. for the distribution exists, with domain of definition having 0 as an interior point, and hence there are obvious advantages.) However, we may stress here that the specific construction that C.R. Rao {\em et al.} (2009) have given to get $M^*$ from $M$ assumes (in our notation)  $v_0= w_0 >0$ (as implied by (4.4) of the cited paper).}
\end{example}
}
\end{remark}

\begin{remark}
\label{rem:3}
{\rm   From the proofs that we have given above, especially for Theorem \ref{th:w-h} and Corollary \ref{col:3}, it is obvious that, for each of Theorem \ref{th:w-h} and Corollaries \ref{col:1}-\ref{col:4}, there exists a sequence $\{ X_{1,n}- X_{2,n}:\, n=1,2,\ldots\}$, with, for each $n$, $X_{1,n}$ and $X_{2,n}$ as independent discrete nonnegative r.v.'s having log-convex distributions (on $\{0, b_n, 2 b_n,\ldots\}$, for some $b_n>0$), converging in distribution to an r.v., with distribution, that is claimed to be i.d., in the respective assertion; from C.R. Rao {\em et al.} (2009, Section 4), it follows further that, in the case of Corollaries \ref{col:2} and \ref{col:4}, the observation remains valid with ``log-convex distributions'' replaced by ``completely monotone distributions or their scale variations''. 
}
\end{remark}

\begin{remark}
\label{rem:4}
{\rm  It follows easily that the classes of the mixtures met in the statements of Corollaries \ref{col:1}, \ref{col:3} and \ref{col:4}, and the class of the distributions in the statement of Corollary \ref{col:2}, are all convex (with members that are i.d). In view of this, it is obvious, amongst other things, that if $X$ is an r.v. with distribution as a member of the class of mixtures in the statement of Corollary \ref{col:3} (or, in particular, in the statement of Corollary \ref{col:4}) and $W$ is an r.v. independent of $X$, then, by the closure property of the class of i.d. distributions, $WX$ and $W |X|$ are i.d. (Moreover, in view of the relevant information in Remark \ref{rem2_2F} $(ii)$ or in the proof of Corollary \ref{col:4}, respectively, Fubini's theorem implies that the distributions of $WX$ and $W|X|$ lie in the class relative to the distribution of $X$, referred to.) It may be worth pointing out in this place that if  $X^*$ is an r.v. with its d.f. as $G$ of Remark \ref{rem:2F}, but with $a=0$, and $W$ is a nonnegative r.v. independent of $X^*$, then, by (the related information in) Remarks \ref{rem2_2F} $(iii)$ and \ref{rem:2F}, $WX^*$ is i.d, since, there is no loss of generality in assuming $W>0$ and we have the d.f. of $WX^*$, in this case, to be a specialized version of $H$ of Remarks \ref{rem2_2F} $(iii)$, with $V=1/W$.}
\end{remark}

\begin{remark}
\label{rem:5}
{\rm In the previous remark, we came across  some cases of independent r.v.'s $W$ and $X$ with $WX$ i.d. One may now raise a question as to whether discrete versions of these relative  to integer-valued r.v.'s are valid. That the answer to this question is in the negative is shown by the following example.

\begin{example}
\label{ex:2}
{\rm Let $W$ and $X$ be independent r.v.'s with $X$ geometric and $W$ a nondegenerate r.v. whose distribution is concentrated on $\{2, 3\}$. Then, since $WX$ is a nonnegative integer-valued r.v. with support of its distribution so that it includes the points 0, 2 and 3, but not the point 5, Theorem 4.2.3 of Steutel (1970) or Theorem II.8.2 of Steutel \& van Harn (2004) implies that this latter r.v. is non-i.d. (In view of the closure property of the class of i.d. distributions, one can, obviously, obtain more general examples to illustrate that this is so.) Incidentally, the present example illustrates also that there exist non-i.d. mixtures of compound geometric distributions.}
\end{example}}
\end{remark}

\begin{remark}
\label{rem:6}
{\rm   Let $X$ and $X^*$ be r.v.'s as in Remark \ref{rem:4}. Also, let $\phi:\RR \mapsto \RR$ and $\phi^*:[0,\infty) \mapsto [0,\infty)$ be one-to-one and onto (i.e., bijective) functions with $\phi(0)=\phi^*(0)=0$, for which their inverses $\phi^{-1}$ and $(\phi^*)^{-1}$ are such that
\begin{equation*}
\label{eq:fanisREm6}
\phi^{-1}(x) = \left\{
\begin{array}{ll}
-\int_0^{|x|} \psi_1(y)dy & {\rm if}\;\;\; x<0,\\
\int_0^x \psi_2(y)dy & {\rm if}\;\;\;x>0, 
\end{array} \right.
\end{equation*}
where 
$\psi_r$,  $r=1,2$,  are decreasing log-convex functions on $(0,\infty)$ with (in standard notation for derivatives) $-\psi_r^{'}(y)$, $y \in (0,\infty)$,  also as log-convex for $r=1,2$, and
$$
(\phi^*)^{-1}(x) = \int_0^x \psi^*(y)dy \quad \text{if} \quad x>0,
$$
where  $\psi^*$ is a decreasing log-convex function on $(0,\infty)$. Then, denoting respectively by $H$ and $H^*$ the d.f.¢s of $\phi(X)$ and $\phi^*(X^*)$, it is seen, for example, that, for $x > 0$, (in standard notation for derivatives)  $- H''(x)$ and $(H^*)'(x)$ are both (nonincreasing) continuous log-convex functions. (To understand this last claim properly, note that if $g$ is proportional to a log-convex density on $(0,\infty)$, then, by Remark \ref{rem2_2F} $(ii)$, $g((\phi)^{-1}(x))$, $x >0$, $g((\phi^*)^{-1}(x))$, $x>0$, and $\int_0^{\infty} g(y+\phi^{-1}(x))\,dy$, $x>0$, are (nonincreasing) continuous log-convex functions, and, also, that the products (by Definition \ref{def1_FF}) and the sums (by Lemma \ref{lem1_FN}) of finitely many log-convex functions on $B$ are log-convex on $B$.) Hence, in view of what we have already pointed out in Remark \ref{rem:4}, it follows in particular (appealing to the underlying symmetry in the case of $H$) that if $W_1$ and $W_2$, with $W_2$ nonnegative, are r.v.'s, independent of $X$ and $X^*$, respectively, then, for all $n\in\{1,2,\ldots\}$, $W_1X^n$ and, for all $\beta \in [1,\infty)$, $W_1|X|^{\beta}$ and $W_2(X^*)^{\beta}$ are i.d. (Clearly, there exist several other possibilities here.)
}
\end{remark}

\begin{remark}
\label{rem:7}
{\rm  If $\{p_x: x=0,1,\ldots\}$ is a probability distribution so that $p_x=u_x-u_{x+1}$, $x=0,1,\ldots$, with $\{u_x: x=0,1,\ldots\}$ completely monotone, then,  by Corollary \ref{col:2}, it is i.d. However, that this assertion does not remain valid, if we replace the condition that $\{u_x\}$ is completely monotone by that it is a {\tt KS}, is shown by the following example.

\begin{example}
\label{ex:3}
{\rm Let $b, c \in (0,1)$ with $b \leq c$, and $\{u_x: x=0,1,\ldots\}$ be so that 
\begin{equation*}
\label{eq:fanisREx3}
u_x= \left\{
\begin{array}{ll}
b^x& {\rm if}\;\;\; x=0,1,\\
b^2c^{x-2} & {\rm if}\;\;\;x=2,3,\ldots~. 
\end{array} \right.
\end{equation*}
Also, let $p_x=u_x-u_{x+1}$, $x=0,1,\ldots$, and $f$ be the corresponding probability generating function (p.g.f.) . Observe now that $\{u_x\}$ considered is a Kaluza sequence and if we choose $c$ sufficiently close to 1 (with $b$ fixed), then, $\{p_x\}$, in this case, turns out to be non-i.d.,  since, in spite of that $0 < f(0) < 1$, we cannot have here  $f(\cdot)= e^{- \lambda + \lambda g(\cdot)}$ with $\lambda > 0$ and $g$ as a p.g.f. satisfying $g(0)=0$. (To see the validity of the claim on $\{p_x\}$, it is sufficient if we verify that $\ln f(s)$ has its second derivative at $s=0$ to be negative.)}
\end{example} }
\end{remark}

\begin{remark}
\label{rem:8}
{\rm  One may now raise a question as to whether Theorem \ref{th:w-h} remains valid if we take in place of (\ref{eq:disc}) its version, obtained from it, replacing the four summations, on its right hand side, respectively, by $v_{-|x|-1}$, $v_{x+1}$, $v_{-1}$ and $v_1$. Taking a hint  from Example \ref{ex:3}, we can now construct the following example to show that the answer to the question is in the negative.

\begin{example}
\label{ex:4}
{\rm Let $\{p_x: x=0,\pm 1, \pm 2, \ldots\}$ be a probability distribution satisfying 
\begin{equation}
\label{eq:ex4FNEW} K p_x = \left\{
\begin{array}{ll}
c^{|x|}  & {\rm if}\;\;\; x=0,-1,-2,\ldots,  \\
b^2 c^{x-2} & {\rm if}\;\;\;
 x=2,3,\ldots, \\
b 
& {\rm if}\;\;\;x=1,
\end{array} \right.
\end{equation}
with $K>0$ and $b, c \in (0,1)$ so that $b \leq c$. Clearly, (\ref{eq:ex4FNEW}) may, now, be viewed as the version of (\ref{eq:disc}), that is sought, with, e.g., $v_1=v_{-1}=1$ and $v_j$'s for $j \neq -1,1$ defined in obvious way. Consequently, if $X$  is an r.v. whose distribution is $\{p_x\}$ then (with $t$ complex)
\begin{eqnarray}
\label{eq:ex4M}
K \EE(e^{tX}) &=& \big(1-ce^{-t}\big)^{-1} +b e^t + b^2 e^{2t} \big(1-ce^t\big)^{-1}, \nonumber\\
&=& (1-bc)\big(1-ce^{-t}\big)^{-1}g(e^t), \quad \text{Re}(t) \in (\ln c, -\ln c),
\end{eqnarray}
with
$g(s) = 1 +bs + b^2s^2 (1-c^2) (1-bc)^{-1}(1-cs)^{-1}$,  for each $s=e^t$, $\text{Re}(t) \in (\ln c, -\ln c)$.
Suppose now that we specialize to the case with $b+b^2(1+c)(1-bc)^{-1}<1$ and $(1-c^2)(1-bc)^{-1}<1/2$
(which subsumes, e.g., the case with $b=1/3$ and $c=8/9$). On applying the standard power-series expansion for $\ln(1+z)$, $|z|<1$, by minor manipulation, (\ref{eq:ex4M}) implies then that, in this special case, there exists a sequence $\{q_x: x=0,\pm 1,\pm2, \ldots\}$ of reals, with $q_0=0$, $q_2 <0$, $q_x=c^{|x|}/|x|$ if $x <0$, and $\sum_{x=2}^{\infty} |q_x| < \infty$, such that
$$
\ln (\EE(e^{tX}))=\sum_{x=-\infty}^{\infty} (e^{tx}-1)q_x, \quad \text{Re}(t) \in (\ln c, 0],
$$
asserting, in view of  Lukacs (1970, Remark 1, p.118) or  Blum \& Rosenblatt (1959, Theorem 1), that $\{p_x\}$ is non-i.d. (Incidentally, Blum \& Rosenblatt (1959, Theorem 1) tells us that any nondegenerate i.d. distribution with at least one discontinuity is so that its characteristic function (ch.f.) has a unique representation as the product of a degenerate ch.f. and a ch.f. that is of the form $e^{-\lambda+\lambda \varphi(\cdot)}$ with $\lambda >0$ and $\varphi$ as the ch.f. of a d.f. $G$ satisfying $G(0)-G(0-)=0$; note that if the i.d. distribution corresponds to an integer-valued r.v., then the degenerate ch.f. and $G$ referred to here also correspond to such r.v.'s.)}
\end{example}
}
\end{remark}

\begin{remark}
\label{rem:9}
{\rm A further example, i.e., Example \ref{ex:5}, that appears below, tells us that Corollary \ref{col:3} does not remain valid if  ``$\int_{|x|}^{\infty}v_1(y)dy$'' and  ``$\int_{x}^{\infty}v_2(y)dy$'' in (\ref{eq:fanisCor3}) are replaced by 
 ``$v_1(|x|)$'' and  ``$v_2(x)$'', respectively. (Incidentally, as a by-product of this, the argument used in the proof of Corollary \ref{col:3}, with obvious alterations, then implies that the answer to the question in Remark \ref{rem:8} is in the negative, thus, supporting the conclusion of Example \ref{ex:4}.)
 
\begin{example} {\rm Let $F$ be an absolutely continuous d.f. with p.d.f. $f$~ such that, for an appropriate constant $K>0$,
\label{ex:5}
\begin{equation}
\label{eq:ex5FNEW} K f(x) = \left\{
\begin{array}{ll}
e^{-\delta |x|}  & {\rm if}\;\;\; x \leq 0,  \\
\alpha e^{-\delta x + h(x)} & {\rm if}\;\;\;
 x > 0,
\end{array} \right.
\end{equation}
with $0< \alpha <e^{-1}$, $\delta>0$ and $h$ as in Example \ref{ex:0F}. If $X$ is an r.v. with d.f.  $F$, we see then, using, in particular, Fubini's theorem or the method of integration by parts, that (with $t$ complex)
\begin{eqnarray}
K \EE(e^{tX}) 
&=& (\delta +t)^{-1} \Big\{ 1+ \alpha(\delta+t) \int_0^{\infty} e^{tx-\delta x + h(x)} dx\Big\}\nonumber \\
&=& (\delta +t)^{-1} \Big\{ 1-\alpha e + \alpha \int_0^{\infty} e^{(\delta+t)x-2\delta x + h(x)} (2\delta-h'(x))dx\Big\}\nonumber \label{eq:lastex}\\
&=& (1-\alpha e)(\delta+t)^{-1} g^*(t),  \quad \text{Re}(t) \in (-\delta,\delta),
\end{eqnarray}
where
$$
g^*(t)= 1+ \alpha (1-\alpha e)^{-1} \int_0^{\infty} e^{tx-\delta x + h(x)} (2\delta-h'(x))dx, \quad \text{Re}(t) \in (-\delta,\delta).
$$
Suppose we now take, for convenience, $\alpha (1-\alpha e)^{-1} < 1/4$ and $\delta=\alpha^2$, and assume that $F$ is i.d. with  L$\acute{\rm e}$vy measure $\nu$. Clearly, in this case, 
$0<g^*(0)-1 < \alpha (1-\alpha e)^{-1} (\delta (e-1) +e+1) <1$,
and the function $\phi$, defined by
$\phi(s) = (g^*(is)-1)/(g^*(0)-1)$, $s \in \RR$,
is a ch.f.; denote by $G$ the d.f. relative to $\phi$. (The inequalities for $g^*(0)-1$, an integral over $(0,\infty)$, are obvious, in view of the assumptions for $\alpha$ and $\delta$, on expressing it as the sum of the appropriate integrals over $(0,1)$ and $[1,\infty)$, respectively.) The table on spectral measures for certain i.d. distributions given by Lukacs (1970, p. 120), and the power-series expansion for $\ln(1+z)$ with $z=(g^*(0)-1)\phi(s)$, $s \in \RR$, in conjunction with the uniqueness theorem for the Fourier transform of a finite signed measure (given, e.g., as Corollary 1.1.2 in Rao \& Shanbhag (1994, p. 2)), imply then, by (\ref{eq:lastex}), that
$$
\nu((-\infty, x)) = \int_{-\infty}^x  |y|^{-1} e^{\delta y} dy, \quad \text{if} \quad x \in (-\infty, 0),
$$
and
\begin{equation}
\label{eq:ex5FNEW2}
\nu((0,x]) = \sum_{n=1}^{\infty} n^{-1}  (-1)^{n-1} (g^*(0)-1)^n G^{n*}(x), \quad \text{if} \quad x \in (0,\infty),
\end{equation}
where, for each $n \in \{1,2,\ldots\}$, $G^{n*}$ is the n-fold convolution of $G$ with itself. From 
(\ref{eq:ex5FNEW2}), we see that there exists a function $o: (0,\infty) \mapsto \RR$ with $\lim_{y \rightarrow 0+}(o(y)/y)=0$, so that, for sufficiently small $\alpha$,
$$
\nu((1,2]) =(1/2)\, \alpha^2 \int_0^1 (e^{h(1-y)} -1) e^{h(y)} h'(y) dy + o(\alpha^2) <0,
$$
leading us to a contradiction, and, hence, supporting the claim made in Remark \ref{rem:9}.}
\end{example}}
\end{remark}

\begin{remark}
\label{rem:10} {\rm To shed further light on the conclusions of Examples \ref{ex:4} and \ref{ex:5}, we may give the following relevant information: Extending $g$ and $g^*$ appearing in these examples appropriately with notation, for convenience, for the extensions respectively as $g_c$ and $g^*_{\delta}$, so that their domains of definition are the sets of complex numbers, respectively, with moduli lying in $(0,1/c)$ and with real parts lying in $(-\infty, \delta)$, it is seen that, for some (real) $t^* \in (-\infty, 0)$, $g_c(\exp\{ i s+ t^*\})/g_c(\exp\{t^*\})$, $s \in \RR$, and $g^*_{\delta}( is + t^*)/g^*_{\delta}(t^*)$, $s \in \RR$, are ch.f.'s. (The extensions referred to here can be assumed to be analytic continuations of their original versions.) If we now allow $c$  and $\delta$ to vary as (distinct) members of a $c$-sequence tending to 1 and a $\delta$-sequence tending to 0, respectively,  then the resulting sequences of ch.f.'s converge to the ch.f.'s of certain nondegenerate bounded r.v.'s.  Clearly, the limiting distributions in the two cases referred to are non-i.d., explaining indirectly, as to why we have the contradictions in the two examples.}
\end{remark}

\section*{\large Acknowledgements}
We are grateful to the two referees for their useful comments on the earlier draft of the paper.


\begin{thebibliography}{99}

\bibitem{}
Blum, J.R. \& Rosenblatt, M. (1959). On the structure of infinitely divisible distributions. {\em Pacific J. Math.}, {\bf 9}, 1--7.

\bibitem{}
Donoghue, W.F., Jr. (1969). {\em Distributions and Fourier Transforms.} New York: Academic Press.

\bibitem{}
Feller, W. (1971). {\em An Introduction to Probability Theory and
its Applications}. Vol. II, 2nd Edition, New York: John Wiley \&
Sons.

\bibitem{}
Goldie, C.M. (1967). A class of infinitely divisible random
variables. {\em Proc.\ Cambridge\ Philos.\ Soc.}, {\bf 63},
1141--1143.

\bibitem{}
Hardy, G., Littlewood, J.E. \& P\'{o}lya, G. (1952). {\em Inequalities}. 2nd Edition, Cambridge University Press, Cambridge. 

\bibitem{}
Kaluza, T. (1928). \"Uber die koeffizienten reziproker potenzreihen. {\em Math.\ Z.},
{\bf 28}, 161--170.

\bibitem{}
Kingman, J.F.C. (1972). {\em Regenerative Phenomena.} New York: John
Wiley \& Sons.

\bibitem{}
Lo$\grave{{\rm e}}$ve, M. (1963). {\em Probability Theory}. 3rd
Edition, Princeton: Van Nostrand.

\bibitem{}
Lukacs, E. (1970). {\em Characteristic Functions}. 2nd Edition, London: Griffin.

\bibitem{}
Rao, C.R. \& Shanbhag, D.N. (1994). {\em Choquet-Deny Type
Functional Equations with Applications to Stochastic Models}.
Chichester: John Wiley \& Sons.

\bibitem{}
Rao, C.R., Shanbhag, D.N., Sapatinas, T. \& Rao, M.B. (2009). Some
properties of extreme stable laws and related infinitely divisible
random variables. {\em J.\ Statist.\ Plann.\ Inference}, {\bf 139},
802--813.

\bibitem{sha}
Shanbhag, D.N. (1977). On renewal sequences. {\em Bull.\ London
Math.\ Soc.}, {\bf 9}, 79--80.

\bibitem{}
Shanbhag, D.N. \& Sreehari, M. (1977). On certain self-decomposable
distributions. {\em Z.\ Wahrsch.\ Verw.\ Gebiete}, {\bf 38},
217--222.

\bibitem{}
Shanbhag, D.N., Pestana, D. \& Sreehari, M. (1977). Some further
results in infinite divisibility. {\em Math.\ Proc.\ Cambridge\
Philos.\ Soc.}, {\bf 82}, 289--295.

\bibitem{}
Steutel, F.W. (1967). Note on the infinite divisibility of
exponential mixtures. {\em Ann.\ Math.\ Statist.}, {\bf 38},
1303--1305.

\bibitem{}
Steutel, F.W. (1970). {\em Preservation of Infinite Divisibility
under Mixing and Related Topics.} Mathematical Centre Tracts, Vol.
{\bf 33}, Amsterdam: Mathematisch Centrum.

\bibitem{}
Steutel, F.W. \& van Harn, K. (2004). {\em Infinite Divisibility of
Probability Distributions on the Real Line}. New York: Marcel
Dekker.

\bibitem{}
Titchmarsh, E.C. (1978). {\em The Theory of Functions}. 2nd Edition, Oxford: Oxford University Press.

\bibitem{}
Zygmund, A. (2002). {\em Trigonometric Series}. Volumes {\bf I} \& {\bf II}, 3rd Edition. Cambridge: Cambridge University Press. 

\end{thebibliography}
\end{document}